\documentclass[11pt]{article}
\usepackage[reqno]{amsmath} 
\usepackage{amssymb,amsthm, enumerate,verbatim,graphicx}

\newcommand\al{\alpha}
\newcommand\be{\beta}

\newcommand\eps{\epsilon}

\newcommand\la{\lambda}

\newcommand\re{{\mathbb R}}

\newcommand\nrm[1]{||#1||}
\newcommand\opk[1]{\mathop{\mathrm{#1}}\nolimits}
\newcommand\one{{\bf1}}
\newcommand\rk{\opk{rk}}

\newcommand\abs[1]{\left|#1\right|}
\newcommand\Rep{\mathrm{Rep}}

\newcommand\defn[1]{\textsl{#1}}
\newtheorem{theorem}{Theorem}
\newtheorem{lemma}[theorem]{Lemma}
\newtheorem{corollary}[theorem]{Corollary}

\theoremstyle{remark}
\newtheorem*{example}{Example}

\begin{document}

\title{Minimal Euclidean representations of graphs}

\author{Aidan Roy\footnote{email:aroy@qis.ucalgary.ca} \\
Department of Mathematics and Statistics \& \\
Institute for Quantum Information Science \\ 
University of Calgary\\ Calgary, Alberta T2N 1N4, Canada}


\maketitle

\begin{abstract}A simple graph $G$ is \defn{representable} in a real vector space of dimension $m$ if there is an embedding of the vertex set in the vector space such that the Euclidean distance between any two distinct vertices is one of only two distinct values $\al$ or $\be$, with distance $\al$ if the vertices are adjacent and distance $\be$ otherwise. The \defn{Euclidean representation number} of $G$ is the smallest dimension in which $G$ is representable. In this note, we bound the Euclidean representation number of a graph using multiplicities of the eigenvalues of the adjacency matrix. We also give an exact formula for the Euclidean representation number using the main angles of the graph.
\end{abstract}

\section{Introduction}\label{sec:intro}

Recently, Nguyen Van Th\'e \cite{NguyenVanThe08} revived a problem of representing graphs in Euclidean space, which, according to Pouzet \cite{Pouzet79}, was originally introduced by Specker before 1972.

A simple graph is \defn{representable} in $\re^m$ if there is an embedding of the vertex set in $\re^m$ and distinct positive constants $\al$ and $\be$ such that for all vertices $u$ and $v$,
\[
\nrm{u-v} = \begin{cases}
\al, & u \sim v; \\
\be, & \text{otherwise.}
\end{cases}
\]
(Here $\nrm{x} := \sqrt{x^Tx}$.) We will call the smallest $m$ such that $G$ is representable in $\re^m$ the \defn{Euclidean representation number} of $G$ and denote it $\Rep(G)$. The complete graph $K_n$ and its complement, the empty graph, are representable in $\re^{n-1}$ via a regular simplex. By a bound on the size of a Euclidean $1$-distance set \cite{Bannai83,Blokhuis84}, they are not representable in smaller dimensions, so $\Rep(K_n) = n-1$. Nguyen Van Th\'e \cite{NguyenVanThe08} showed that if $G$ is a graph on $n$ vertices that is not the complete graph or empty graph, then $\Rep(G) \leq n-2$. Here, we use the multiplicities of the smallest and second smallest eigenvalues of the adjacency matrix $A(G)$ to give upper and lower bounds for the representation number. The main result is the following.

\begin{theorem}
\label{thm:summary}
Let be $G$ be a graph on $n$ vertices which is not complete or empty. If $G$ or its complement is the disjoint union of complete graphs, $r$ of which are of maximal size, then 
\[
\Rep(G) = n-\max\{r,2\}.
\]
Otherwise, let $m_1$ and $m_2$ respectively denote the multiplicity of the smallest and second smallest eigenvalue of the adjacency matrix of $G$, and similarly define $\overline{m_1}$ and $\overline{m_2}$ for the complement of $G$. Then
\[
n-1-\max\{m_1,\overline{m_1},m_2+1,\overline{m_2}+1\} \; \leq \; \Rep(G) \; \leq \; n-\max\{m_1,\overline{m_1},2\}.
\]
If $G$ is regular, then
\[
\Rep(G) = n-1-\max\{m_1,\overline{m_1}\}.
\]
\end{theorem}

To get a precise characterization for irregular graphs, we must consider not just the eigenvalues but the actual eigenspaces, or at the very least, the main angles of eigenspaces. An exact formula for the representation number is given below in Theorem \ref{thm:representation}. 

It is not surprising that the Euclidean representation number of a graph is closely related to the multiplicity of the smallest eigenvalue. If $\tau_1$ is the smallest eigenvalue of $A(G)$, then the positive semidefinite matrix $A(G) - \tau_1 I$ is the Gram matrix of a set of $n$ vectors in $\re^{n-m_1}$. This technique is perhaps the most common method of embedding a graph in a vector space and it plays a critical role in, for example, the characterization of graphs with least eigenvalue $\tau_1 \geq -2$ and the theory of two-graphs and equiangular lines \cite{Cameron91,Cvetkovic95,Godsil01}.

Any representation of a graph in $\re^m$ is by definition a \defn{Euclidean $2$-distance set}: a set of vectors such that only two nontrivial distances occur between vectors \cite{Bannai83,Blokhuis84}. Conversely, the distance relation of any Euclidean $2$-distance set of size $n$ defines an $n$-vertex graph. The problem of finding the largest $2$-distance set (and its corresponding graph) in a given dimension $m$ has been well studied, but to our knowledge the issue in this paper --- finding the smallest dimension for a given graph --- has not previously been resolved. On the other hand, the main tool used in this paper, namely the Euclidean distance matrix, has been applied to $s$-distance sets quite often (see for example Larman, Rogers, and Seidel \cite{Larman77}). Moreover, there is a related problem in graph network theory that has received considerable attention: given a graph, and specifying the length of each edge, find an embedding of the graph in $\re^m$ \cite{Alfakih98}. That problem is NP-hard \cite{Saxe79}.


Theorem \ref{thm:summary} says in particular that a graph with all simple eigenvalues (eigenvalues with multiplicity $1$) has Euclidean representation number $n-2$ or $n-3$. Since eigenvalues of a random graph are almost surely simple, Nguyen Van Th\'e's upper bound of $\Rep(G) \leq n-2$ is often tight. (For the distribution of eigenvalues in random graphs see Chung, Lu \& Vu \cite{Chung03}.) In fact, Theorem \ref{thm:representation} shows that almost surely, the representation number of a random graph is $n-2$. It seems that only special graphs, such as line graphs or graphs with a high degree of regularity or symmetry, can be represented in smaller dimensions. 

\section{General graphs}\label{sec:general}

We characterize representations of $G$ using Euclidean distance matrices. An $n \times n$ matrix $M$ is a \defn{Euclidean distance matrix} of a set of $n$ vectors in $\re^m$ if the rows and columns of $M$ are indexed by the vectors, and $M_{u,v} = ||u-v||^2$. Thus every Euclidean distance matrix is nonnegative with zero diagonal.

Let $A := A(G)$ denote the adjacency matrix of $G$, let $J$ denote the all-ones matrix and let $\overline{A} := J-I-A$ be the adjacency matrix of the complement of $G$. Then up to scaling, the Euclidean distance matrix of a representation of $G$ has the form $M = A + b\overline{A}$, for some $b > 0$ and $b \neq 1$. The smallest dimension $m$ such that $M$ is the distance matrix of a set of points in $\re^m$ is called the \defn{embedding dimension} of $M$ or \defn{dimensionality} of $M$ \cite{Gower85,Hayden91}. Therefore the Euclidean representation number of $G$ is the smallest embedding dimension of a Euclidean distance matrix representing $G$.

We will use two characterizations of Euclidean distance matrices: the first is due to Schoenberg \cite{Schoenberg38}. Let $\one$ denote the all-ones vector and let $P := I - \frac{1}{n}\one\one^T$, the projection matrix for the space $\one^\perp$.

\begin{theorem}
\label{thm:edmchar}
Let $M$ be a symmetric matrix zero diagonal and positive off-diagonal entries. Then $M$ is a Euclidean distance matrix if and only if $PMP$ is negative semidefinite (that is, $M$ is negative semidefinite on $\one ^\perp$). The embedding dimension of $M$ is the rank of $PMP$.
\end{theorem}

The second characterization is a generalization of Schoenberg's due to Gower \cite{Gower85,Gower82}.

\begin{theorem}
Let $M$ be a symmetric matrix zero diagonal and positive off-diagonal entries. For any real vector $v$ such that $v^T\one = 1$, let 
\[
F := (I - \one v^T)M(I - v\one^T).
\]
Then $M$ is a Euclidean distance matrix if and only if $F$ is negative semidefinite. The embedding dimension of $M$ is the rank of $F$.
\label{thm:edmchar2}
\end{theorem}


In order to give a complete description of the representation number, we will also use the main angles \cite[Chapter 5]{Cvetkovic97} of a graph $G$. It will be convenient to order the distinct eigenvalues of $A(G)$ from smallest to largest as $\tau_1,\tau_2,\ldots,\tau_s$. Given an eigenvalue $\tau_i$, with an eigenspace $E_i$ and projection matrix $P_i$ onto that eigenspace, the \defn{main angle} of $\tau_i$ is 
\[
\be_i := \frac{1}{\sqrt{n}}\nrm{P_i \one}.
\]
Note that $0 \leq \be_i \leq 1$, and $\be_i = 0$ if and only if $E_i \subseteq \one^\perp$. In particular $\be_s > 0$, as the largest eigenvalue $\tau_s$ has an eigenvector with nonnegative entries by the Perron-Frobenius Theorem (see \cite[Theorem 0.2]{Cvetkovic95} or \cite[Theorem 8.2.11]{Horn85}). Also note that if $G$ has $n$ vertices, then
\[
n \be_i^2 = \max_{\substack{v \in E_i \\ v^Tv = 1}} (v^T\one)^2 = \max_{v \in E_i} \frac{(v^T\one)^2}{v^Tv}
\]
and likewise
\begin{equation}
\frac{1}{n \be_i^2} = \min_{v \in E_i} \frac{v^Tv}{(v^T\one)^2} = \min_{\substack{v \in E_i \\ v^T\one = 1}} v^Tv.
\label{eqn:mainvaluechar}
\end{equation}
If $\be_i \neq 0$, then $\tau_i$ is called a \defn{main eigenvalue} of $G$. 

We now consider Euclidean distance matrices of the form $M = A + b\overline{A}$. If $b = 1$, then $M = J-I$, which is not the distance matrix of a valid representation of a graph unless $G$ is the complete or empty graph. Therefore there are two separate cases: $b > 1$ and $0< b < 1$. Without loss of generality, we may assume that $b > 1$, provided that we also consider Euclidean representations of the complement of $G$.

Any graph which is not a disjoint union of complete graphs has at least one eigenvalue smaller than $-1$ (see \cite[Exercise 3.1]{Cameron91} or \cite[Chapter 1, Proposition 6.1]{Beineke04}).

\begin{lemma}
\label{lem:edmtau}
Let $G$ be a graph which is not the disjoint union of complete graphs and let $\tau_1 < -1$ be the smallest eigenvalue of $A$, with multiplicity $m_1$ and main value $\be_1$. If $b = \tau_1/(\tau_1+1)$, then $M = A + b\overline{A}$ is a Euclidean distance matrix. The embedding dimension of $M$ is $n-m_1-1$ if $\be_1 = 0$ and $n-m_1$ otherwise.
\end{lemma}

\begin{proof}
If $b = \tau_1/(\tau_1+1)$ and $\tau_1 < -1$, then $b > 1$. Let
\[
M = A + b\overline{A} = (1-b)A - bI + bJ.
\]
Since $b > 0$, $M$ has positive off-diagonal entries. Therefore $M$ is a Euclidean distance matrix if and only if $PMP$ is negative semidefinite. 
Instead of $M$, consider
\[
X := (1-b)A - bI.
\]
$X$ has largest eigenvalue $(1-b)\tau_1 - b = 0$. Therefore $X$ is negative semidefinite, so $PXP = PMP$ is also negative semidefinite. Thus $M$ is a Euclidean distance matrix. 

By Theorem \ref{thm:edmchar}, the embedding dimension of $M$ is the rank of $PMP = PXP$. Since $X$ is negative semidefinite, say $X = -U^TU$, the rank of $PXP = -PU^TUP$ is the rank of $UP$. The null space of $U$ is $N(X) = E_1$, the eigenspace of $\tau_1$. Therefore the null space of $UP$ is $\langle\one\rangle \oplus (E_1 \cap \one^\perp)$. The dimension of this space is $m_1+1$ if $\be_1 = 0$ and $m_1$ otherwise.
\end{proof}

\begin{example}
Consider $G = \overline{C_6}$, the complement of the cycle graph on six vertices. The eigenvalues of this graph are $\{3, 1, 0, 0, -2, -2\}$, so $m_1 = 2$. Moreover, $\one$ is an eigenvector for the eigenvalue $3$, so $E_1$ is contained in $\one^\perp$ and $\be_1 = 0$. By Lemma \ref{lem:edmtau}, $M = A + 2\overline{A}$ is a Euclidean distance matrix of embedding dimension $3$. On the other hand, consider $G = 2P_3$, the disjoint union of 2 paths of 3 vertices each. The eigenvalues are $\{\sqrt{2},\sqrt{2},0,0,-\sqrt{2},-\sqrt{2}\}$, so again $m_1 = 2$. However, $E_1$ is not contained in $\one^\perp$, so $\be_1 > 0$ and the embedding dimension of $M = A + (2+\sqrt{2})\overline{A}$ is $4$. 
\end{example}


In addition to $b = \tau_1/(\tau_1+1)$, there is a second choice of $b$ that sometimes results in a distance matrix with small embedding dimension: $b = \tau_2/(\tau_2+1)$.

\begin{lemma}
\label{lem:edmtau2}
Let $\tau_i$ denote the $i$-th smallest distinct eigenvalue of $A$, with multiplicity $m_i$ and main angle $\be_i$. Assume $\tau_2 < -1$, and let $b = \tau_2/(\tau_2+1)$ and $M = A + b\overline{A}$. Then $M$ is a Euclidean distance matrix if and only if $m_1 = 1$, $\be_2 = 0$, and
\begin{equation}
\frac{\be_1^2}{\tau_2-\tau_1} \geq \sum_{i \geq 3} \frac{\be_i^2}{\tau_i-\tau_2}.
\label{eqn:tau2cond}
\end{equation}
Moreover, the embedding dimension of $M$ is either $n-m_2-2$ if equality holds in \eqref{eqn:tau2cond} and $n-m_2-1$ otherwise.
\end{lemma}

\begin{proof} 
Let $\eps_i := \tau_i(1-b) - b = (\tau_i - \tau_2)/(\tau_2+1)$. Since $\tau_i$ is the $i$-th smallest eigenvalue of $A$, $\eps_i$ is the $i$-th largest eigenvalue of $X = (1-b)A - bI$. Let $E_i$ denote the eigenspace of $\eps_i$. 

Since $\eps_1 > 0$ and $\eps_2 = 0$, we have $x^TXx \geq 0$ for every $x$ in $E_1 \oplus E_2$, with equality if and only if $x$ is in $E_2$. If $M$ is a Euclidean distance matrix, then $X$ is negative semidefinite on $\one^\perp$. So, any $x \in (E_1 \oplus E_2) \cap \one^\perp$ satisfies $x^TXx = 0$ and is therefore in $E_2$. That is, $(E_1 \oplus E_2) \cap \one^\perp \subseteq E_2$. Since the dimension of $(E_1 \oplus E_2) \cap \one^\perp$ is at least $\dim(E_2) + \dim(E_1) - 1 \geq \dim(E_2)$, it follows that $(E_1 \oplus E_2) \cap \one^\perp = E_2$. Therefore $\dim(E_1) = 1$, $E_1 \nsubseteq \one^\perp$, and $E_2 \subseteq \one^\perp$. In other words, $m_1 = 1$, $\be_1 \neq 0$, and $\be_2 = 0$.

Now, suppose $m_1 = 1$, $\be_1 \neq 0$, and $\be_2 = 0$. Let $v_1$ be the eigenvalue in $E_1$, normalized so that $v_1^T\one = 1$. By Theorem \ref{thm:edmchar2}, $M$ is a Euclidean distance matrix if and only if
\[
F := (I - \one v_1^T)M(I - v_1\one^T) = (I - \one v_1^T)X(I - v_1\one^T)
\]
is negative semidefinite. Note that $Fv_1 = 0$, and for any $v_2 \in E_2 \subseteq \one^\perp$, $Fv_2 = 0$. Therefore is suffices to consider if $x^TFx \leq 0$ for $x$ of the form $x = \sum_{i\geq3} a_iv_i$, with $v_i \in E_i$. Without loss of generality, assume $x^T\one = 1$ and $v_i^T \one = 1$, so $\sum_{i \geq 3} a_i = 1$. Then 
\[
x^TFx = \sum_{i \geq 3}a_i^2\eps_i(v_i^Tv_i) + \eps_1(v_1^Tv_1).
\]
Using the main values of $A$ as in \eqref{eqn:mainvaluechar}, the smallest value of $v_i^Tv_i$ ($i \geq 3$) is $1/n\be_i^2$. So, $F$ is negative semidefinite if and only if
\begin{equation}
\max_a \sum_{i \geq 3}a_i^2\eps_i/\be_i^2 + \eps_1/\be_1^2 \text{ subject to } \sum_{i \geq 3} a_i = 1
\label{eqn:quadmax}
\end{equation}
is nonpositive. This quadratic optimization problem (see for example \cite[Chapter 9]{Wismer78}) obtains its maximum at $a_i = -c\be_i^2/\eps_i$, where $c$ is a normalization constant $c = 1/(\sum_{i \geq 3} -\be_i^2/\eps_i)^2$. Plugging this maximum into \eqref{eqn:quadmax}, we find that $M$ is a Euclidean distance matrix if and only if
\[
\frac{1}{-\sum_{i \geq 3}\be_i^2/\eps_i} \geq \eps_1/\be_1^2,
\]
from which \eqref{eqn:tau2cond} follows.

If equality holds in \eqref{eqn:tau2cond}, then the choice of $x = \sum_{i\geq3} a_iv_i$ with equality satisfies $x^TFx = 0$ and therefore $Fx = 0$. In this case, the null space of $F$ contains $x$, $E_1$, and $E_2$, so the embedding dimension is $n-m_2-2$. Otherwise, the null space of $F$ is $E_1 \oplus E_2$ and the embedding dimension is $n-m_2-1$.
\end{proof}
 
\begin{figure}
\center{\includegraphics[width=100mm]{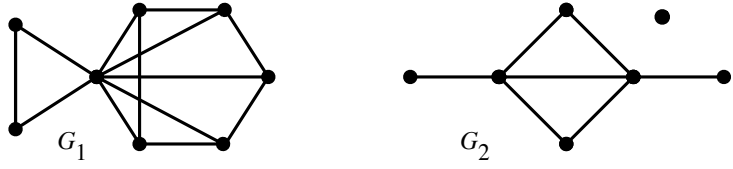}}
\caption{graphs $G_1$ and $G_2$ have representations in $\re^{n-m_2-1}$ and $\re^{n-m_2-2}$ respectively.\label{fig:sevenverts}}
\end{figure}

\begin{example}
Consider $G = G_1$ in Figure \ref{fig:sevenverts}. The eigenvalues of $G$ are $\{3.694$, $1.252,0.618, -1, -1.618, -1.946\}$, with multiplicities $\{1, 1, 2, 1, 2, 1\}$ and main values $\{0.955, 0.124, 0, 0, 0, 0.269\}$ respectively. It follows that inequality \eqref{eqn:tau2cond} holds without equality, so by Lemma \ref{lem:edmtau2}, $M = A + 2.618\overline{A}$ is a Euclidean distance matrix of embedding dimension $n-m_2-1 = 5$. On the other hand, consider $G = G_2$ in Figure \ref{fig:sevenverts}. The eigenvalues are $\{(1+\sqrt{21})/2$, $(-1+\sqrt{5})/2, 0, (-1-\sqrt{5})/2, (1-\sqrt{21})/2\}$, with multiplicities $\{1,1, 3, 1, 1\}$ and main values 
\[
\Big\{\sqrt{\tfrac{10 + 2\sqrt{21}}{21 + \sqrt{21}}}, 0, \sqrt{\tfrac{1}{5}}, 0, \sqrt{\tfrac{34 - 6\sqrt{21}}{105 + 5\sqrt{21}}}\Big\}
\]
respectively. Equality holds in \eqref{eqn:tau2cond}, so $M = A + \tfrac{1}{2}(3+\sqrt{5})\overline{A}$ is a Euclidean distance matrix of embedding dimension $n-m_2-2 = 4$.
\end{example}

Next, we show that a Euclidean distance matrix of the form $M = A + b\overline{A}$ can not have small embedding dimension unless $b = \tau_1/(\tau_1+1)$ or $b = \tau_2/(\tau_2+1)$. (For related results about eigenvalues of Euclidean distance matrices, see \cite{Larman77,Neumaier81}.)

\begin{lemma}
\label{lem:edmsmaller}
Let $b > 1$ and let $M = A + b\overline{A}$ be a Euclidean distance matrix. If $\tau_2$ denotes the second smallest distinct eigenvalue of $A$ and $\tau_2 < -1$, then $b \leq \tau_2/(\tau_2+1)$. Moreover, the embedding dimension of $M$ is at least $\rk(X)-2$, where $X = (1-b)A - bI$.
\end{lemma}

\begin{proof} 
Let $\eps_1 := \tau_1(1-b) - b$ and $\eps_2 := \tau_2(1-b) - b$. Since $\tau_1$ and $\tau_2$ are the smallest and second smallest eigenvalues of $A$, $\eps_1$ and $\eps_2$ are largest and second largest eigenvalues of $X = (1-b)A - bI$. Now $M$ is a Euclidean distance matrix, so $X$ is negative semidefinite on $\one^\perp$. Since $\one^\perp$ is a space of dimension $n-1$, it follows that $X$ can have at most one positive eigenvalue. Thus $\eps_2 \leq 0$, which implies $b \leq \tau_2/(\tau_2+1)$.

The embedding dimension of $M$ is the rank of $PMP = PXP$. Since $P$ has rank $n-1$, we see that $\rk(XP) \geq \rk(X) - 1$ and $\rk(PXP) \geq \rk(XP)-1 \geq \rk(X) -2$.
\end{proof}

We combine the results of Lemmas \ref{lem:edmtau}, \ref{lem:edmtau2} and \ref{lem:edmsmaller} to find the choice of $b>1$ in $M = A + b\overline{A}$ that gives the best representation of $G$. Nguyen Van Th\'e \cite{NguyenVanThe08} showed that $\Rep(G) \leq n-2$. Lemma \ref{lem:edmsmaller} shows that in order of find a smaller representation, we must choose $b$ such that $X = (1-b)A - bI$ has less than full rank. But $\rk(X) < n$ only if $b = \tau_i/(\tau_i+1)$, where $\tau_i$ is some eigenvalue of $A$. Moreover, by Lemma \ref{lem:edmsmaller}, $M$ is a Euclidean distance matrix with $b > 1$ only if $b \leq \tau_2/(\tau_2+1)$. So the only possible choices of $b$ are $\tau_1/(\tau_1+1)$ and $\tau_2/(\tau_2+1)$. When $b = \tau_1/(\tau_1+1)$, Lemma \ref{lem:edmtau} shows that $M$ is a Euclidean distance matrix and the embedding dimension is either $n-m_1$ or $n-m_1-1$. When $b = \tau_2/(\tau_2+1)$, Lemma \ref{lem:edmtau2} shows that $M$ is sometimes a Euclidean distance matrix, with embedding dimension of $M$ is either $n-m_2-1$ or $n-m_2-2$.

By considering both $G$ and its complement, this gives a complete description of Euclidean representation number of $G$, which we now summarize as a theorem.

Recall that for any graph, the smallest eigenvalue $\tau_1$ is at most $-1$, with equality if and only if $G$ is a disjoint union of complete graphs. Now consider the largest component of such a disjoint union (that is, the complete subgraph with the most number of vertices). If $G$ has $r$ such components of largest size, then the largest eigenvalue of $A$ has multiplicity $r$ and the smallest eigenvalue $\overline{\tau_1}$ of $\overline{A}$ has multiplicity $r-1$. Moreover, it is not difficult to verify that the eigenspace of $\overline{\tau_1}$ is orthogonal to $\one$, so the Euclidean distance matrix given in Lemma \ref{lem:edmtau} for the complement of $G$ has embedding dimension $n-(r-1)-1 = n-r$. This representation is optimal provided that $r \geq 2$. When $G$ is not a disjoint union of complete graphs, the optimal representation either has dimension $n-2$ or occurs when $b \in \{\tau_1/(\tau_1+1)$, $\tau_2/(\tau_2+1)\}$, for either $G$ or its complement.

\begin{theorem}
\label{thm:representation}
Let $G$ be a graph on $n$ vertices and $e$ edges. If $G$ or its complement is the complete graph, then
\[
\Rep(G) = n-1.
\]
If $G$ or its complement is the disjoint union of at least $2$ complete graphs, $r$ of which are of largest size, then
\[
\Rep(G) = n-\max\{r,2\}.
\]
Otherwise, let $\tau_i$ be the $i$-th smallest distinct eigenvalue of $A = A(G)$, with multiplicity $m_i$ and main value $\be_i$. Define
\[
m_1' := \begin{cases}
m_1+1, & \text{if } \be_1 = 0; \\
m_1, & \mbox{otherwise,} 
\end{cases}
\]
and
\[
m_2' := \begin{cases}
m_2 + 2, & \text{if } \tau_2 < -1, m_1 = 1, \be_2 = 0, \text{ and } \tfrac{\be_1^2}{\tau_2-\tau_1} = \sum_{i \geq 3} \tfrac{\be_i^2}{\tau_i-\tau_2}; \\
m_2 + 1, & \text{if } \tau_2 < -1, m_1 = 1, \be_2 = 0, \text{ and } \tfrac{\be_1^2}{\tau_2-\tau_1} > \sum_{i \geq 3} \tfrac{\be_i^2}{\tau_i-\tau_2}; \\
0, & \mbox{otherwise.} 
\end{cases}
\]
Similarly define $\overline{m_1}'$ and $\overline{m_2}'$ for the complement of $G$. Then
\[
\Rep(G) = n-\max\{m_1',m_2',\overline{m_1}',\overline{m_2}',2\}.
\]
\end{theorem}

As an application of Theorem \ref{thm:representation}, consider line graphs. If $G$ is a graph with $n$ vertices and $e$ edges, then the \defn{line graph} $L(G)$ is the graph whose vertices are the edges of $G$, with two edges adjacent if and only if they share a common vertex in $G$. The line graph has $e$ vertices and has smallest eigenvalue $\tau_1 = -2$, with multiplicity $m_1$ at least $e-n$. (This result is due to Sachs \cite{Sachs67}; see also \cite[Theorem 3.8]{Biggs93}.) If follows that $L(G)$ is representable in $\re^n$. More precisely, let $B$ denote the $n \times e$ unoriented incidence matrix of $G$ (the matrix in which $B_{ij} = 1$ if vertex $i$ is incident with edge $j$ and $B_{ij} = 0$ otherwise). Then the adjacency matrix of $L(G)$ is $B^TB - 2I$, so $m_1 = e-\rk(B)$. Moreover, the rank of the incidence matrix of $G$ is easy compute: if $r$ denotes the number of connected components of $G$ that are bipartite, then the rank of $B$ is $n-r$ (see \cite[Theorem 8.2.1]{Godsil01}). Finally, the $\tau_1$-eigenspace $E_1$ of $L(G)$ is orthogonal to $\one$ \cite[Corollary to Theorem 3.38]{Cvetkovic95}, so the representation number of $G$ is $e-m_1-1$, provided that $m_1 \geq \overline{m_1}$. We have:

\begin{corollary}
Let $G$ be a graph with $n$ vertices, $e$ edges, and $r$ bipartite connected components. Then
\[
\Rep(L(G)) \leq n-1-r.
\]
If $e \geq 2(n-r)$, then $\Rep(L(G)) = n-1-r$.
\end{corollary}

Theorem \ref{thm:representation} gives the Euclidean representation number of $G$ in terms of the eigenspaces of $G$ and its complement. In general, there is a relationship between the eigenspace of $\tau_i$ for $A$ and the eigenspace of $-\tau_i-1$ for $\overline{A}$: their dimensions differ by at most one (see \cite{Cvetkovic71} or \cite[Theorem 2.5]{Cvetkovic95}). More generally, the spectrum of the complement of $G$ is determined by the eigenvalues, multiplicities, and main angles of $G$. The following observation is due to Cvetkovic and Doob (\cite{Cvetkovic85}, see also \cite[Proposition 4.5.2]{Cvetkovic97}).

\begin{theorem}
\label{thm:complement}
Let $G$ be a graph with distinct eigenvalues $\{\tau_i\}$, multiplicities $m_i$, and main angles $\be_i$, so the characteristic polynomial of $G$ is $P_G(x) = \prod_i (x-\tau_i)^{m_i}$. Then
\[
P_{\overline{G}}(x) = (-1)^nP_G(-x-1) \left(1 - n\sum_i \frac{\be_i^2}{x+1+\tau_i}\right).
\]
\end{theorem}

\section{Regular graphs} \label{sec:regular}

A graph is \defn{regular} if every vertex has the same degree (number of neighbours). For regular graphs, the formula for the Euclidean representation number of a graph in Theorem \ref{thm:representation} depends only on the multiplicities of the eigenvalues rather than both the multiplicities and the main angles.

If $G$ is regular and the disjoint union of $r$ complete graphs, then $G = rK_{n/r}$. Here $G$ and its complement have representation number $n-r$ (each $K_{n/r}$ is represented in $\re^{n/r-1}$ as a regular simplex, and each component is embedded into a distinct orthogonal subspace of $\re^{n-r}$.) Otherwise, assume $G$ is regular of degree $k$. Then $A\one = k\one$, so $\one$ is an eigenvector of $A$. It follows that the eigenspace of $\tau_1$ is contained in $\one^\perp$ and $\be_1 = 0$. From Theorem \ref{thm:representation}, we see that the representation number any regular $G$ which is not a disjoint union of complete graphs or its complement is 
\[
\Rep(G) = n-1-\max\{m_1,\overline{m_1}\}.
\]

Moreover, the eigenvalues of $A$ precisely determine the eigenvalues of $\overline{A}$. The following now standard result is due to Sachs \cite{Sachs62}. (The result follows from Lemma \ref{thm:complement}, but for a combinatorial proof, see \cite[Theorem 2.6]{Cvetkovic95}).

\begin{theorem}
Let $G$ be an $n$-vertex $k$-regular graph with (not necessarily distinct) eigenvalues $\la_1 = k, \la_2,\ldots,\la_n$, where $\la_i \geq \la_{i+1}$. Then the eigenvalues of $\overline{A}$ are $n-k-1$,$-\la_n-1$,\ldots,$-\la_2-1$.
\end{theorem}


As before, denote the $s$ distinct eigenvalues of $G$ by $\tau_1 < \ldots < \tau_s = k$ with multiplicities $m_1,\ldots,m_s$ respectively. If $G$ is connected, then the multiplicity $\overline{m_1}$ of the smallest eigenvalue of $\overline{A}$ is exactly $m_{s-1}$, the multiplicity of the second largest eigenvalue of $A$. If $G$ is disconnected, then $\overline{m_1} = m_s-1$, where $m_s$ is the multiplicity of $\tau_s = k$ and is also the number of components in $G$. Combining these observations, we have:

\begin{theorem}
\label{thm:regrepresentation}
Let be $G$ be a regular graph on $n$ vertices. If $G = rK_{n/r}$ ($r \neq n$) or its complement, then
\[
\Rep(G) = n-r.
\]
Otherwise, let $m_1$ and $m_{s-1}$ be the multiplicities of the smallest and second largest distinct eigenvalues of $G$. If $G$ is connected, then
\[
\Rep(G) = n-1-\max\{m_1,m_{s-1}\}.
\]
If $G$ is disconnected with $r$ components, then 
\[
\Rep(G) = n-1-\max\{m_1,r-1\}.
\]
\end{theorem}

For more information about the smallest and second largest eigenvalues, see the survey by Seidel \cite{Seidel89}. 

By way of application, consider strongly regular graphs; for background, see \cite{Cvetkovic95}, \cite{Godsil01}, or \cite{Biggs93}. A strongly regular graph $G$ with parameters $(n,k,\la,\mu)$ is a graph with $n$ vertices and valency $k$ such that two distinct vertices have $\la$ common neighbours if they are adjacent and $\mu$ common neighbours otherwise. Such a graph has only two eigenvalues other than $k$, which have multiplicities
\[
m_1,m_2 = \frac{1}{2} \left(n-1 \pm \frac{(n-1)(\mu-\la) - 2k}{\sqrt{(\mu-\la)^2 + 4(k-\mu)}}\right).
\]
The only disconnected strongly regular graph is $rK_{n/r}$, which has representation number $n-r$. Excluding $rK_{n/r}$ and its complement, we have the following:

\begin{corollary}
Let $G$ be a strongly-regular graph which is not a complete multipartite graph or its complement, and has parameters $(n,k,\la,\mu)$. Then
\[
\Rep(G) = \frac{1}{2} \left(n-1 - \frac{\abs{(n-1)(\mu-\la) - 2k}}{\sqrt{(\mu-\la)^2 + 4(k-\mu)}}\right).
\]
\end{corollary}

\begin{example}
The Petersen graph has eigenvalues $3$, $1$, and $-2$, with multiplicities $1$, $5$, and $4$ respectively. Therefore it is representable in $\re^4$ \cite{Lisonek97}.
\end{example}

More generally, given the intersection array of a distance-regular graph $G$, one may readily compute the eigenvalues and multiplicities of $G$ and use Theorem \ref{thm:regrepresentation} to find the Euclidean representation number of $G$. For example, $C_n$, the cycle on $n$ vertices, is a distance-regular graph whose second largest eigenvalue has multiplicity $2$. So for $n \geq 5$, $C_n$ is representable in $\re^{n-3}$.

As another application, consider the following theorem of Petersdorf and Sachs (\cite{Petersdorf70}, see also \cite[Proposition 16.6]{Biggs93}). A graph is \defn{vertex-transitive} is its automorphism group acts transitively on the vertices.

\begin{theorem}
\label{thm:vertextransitive}
Let $G$ be a vertex-transitive graph of degree $k$ and $n$ vertices, and let $\tau_i$ be a simple eigenvalue of $G$. If $n$ is odd, then $\tau_i = k$. If $n$ is even, then $\tau_i = 2\al - k$, where $\al$ is an integer between $0$ and $k$. 
\end{theorem}

\begin{corollary}
If $G$ is a vertex-transitive graph other than the complete or empty graph, and $G$ has an odd number of vertices $n$, then $G$ is representable in $\re^{n-3}$.
\end{corollary}

\begin{proof}
The smallest eigenvalue of $G$ has multiplicity at least $2$. If $G$ is $rK_{n/r}$ or its complement, then $r$ is an odd number, so $r \geq 3$ and therefore $\Rep(G) = n-r \leq n-3$. Otherwise, the result follows from Theorems \ref{thm:regrepresentation} and \ref{thm:vertextransitive}.
\end{proof}

In a similar vein, a graph $G$ is \defn{symmetric} if for all vertices $w,x,y,z$, if $w \sim x$ and $y \sim z$, then there is a an automorphism of $G$ mapping $w$ to $y$ and $x$ to $z$. Every symmetric graph is vertex-transitive. The following theorem due to Biggs \cite[Proposition 16.7]{Biggs93}.

\begin{theorem}
If $G$ is a symmetric graph of degree $k$ and $\tau_i$ is a simple eigenvalue of $G$, then $\tau_i = \pm k$.
\end{theorem}

The eigenvalue $\tau_1 = -k$ occurs in a $k$-regular graph if and only if the graph is bipartite.

\begin{corollary}
If $G$ is a non-bipartite symmetric graph on $n$ vertices, and $G$ is not $K_n$ or $2K_{n/2}$, then $G$ is representable in $\re^{n-3}$. 
\end{corollary}


\section{Acknowledgements}
The author would like to thank Lionel Nguyen Van Th\'e, Chris Godsil, and an anonymous referee for their helpful input. This research is supported by NSERC, MITACS, PIMS, iCORE, and the University of Calgary Department of Mathematics and Statistics Postdoctoral Program. 


\end{document}